\documentclass[10pt]{elsarticle}
\usepackage[cp1251]{inputenc}
\usepackage[english]{babel}
\usepackage{amsmath}
\usepackage{amssymb}
\usepackage{amsfonts}
\usepackage{rotating}
\usepackage{graphicx}
\usepackage{floatflt,epsfig}
\usepackage{lineno,hyperref}
\usepackage{enumerate}
\usepackage{colortbl}
\usepackage[table]{xcolor}
\usepackage{array,tabularx,tabulary,booktabs,ytableau}
\usepackage{longtable}
\usepackage{multirow}
\usepackage{wrapfig}
\usepackage{subcaption}

\newcolumntype{$}{>{\global\let\currentrowstyle\relax}}
\newcolumntype{^}{>{\currentrowstyle}}

\journal{arXiv}
\newtheorem{lemma}{Lemma}

\newtheorem{theorem}{Theorem}

\newtheorem{proposition}{Proposition}
\newcommand{\proof}{\medskip\noindent{\bf Proof.~}}
\begin{document}
\renewcommand{\abstractname}{Abstract}
\renewcommand{\refname}{References}
\renewcommand{\tablename}{Figure.}
\renewcommand{\arraystretch}{0.9}
\thispagestyle{empty}
\sloppy

\begin{frontmatter}
\title{$PI$-eigenfunctions of the Star graphs\tnoteref{grant}}
\tnotetext[grant]{The reported study was funded by RFBR according to the research project 17-51-560008.
The first author is partially supported by the NSFC (11671258) and STCSM (17690740800). The first and the fourth authors
are partially supported by RFBR according to the research project 16-31-00316.
}

\author[02,03,05]{Sergey~Goryainov\corref{cor1}}
\cortext[cor1]{Corresponding author}
\ead{44g@mail.ru}

\author[03]{Vladislav~Kabanov}
\ead{vvk@imm.uran.ru}

\author[01,04]{Elena~Konstantinova}
\ead{e\_konsta@math.nsc.ru}

\author[03,05]{Leonid~Shalaginov}
\ead{44sh@mail.ru}

\author[01]{Alexandr~Valyuzhenich}
\ead{graphkiper@mail.ru}

\address[01]{Sobolev Institute of Mathematics, Ak. Koptyug av. 4, Novosibirsk, 630090, Russia}
\address[02]{Shanghai Jiao Tong University, 800 Dongchuan RD. Minhang District, Shanghai, China}
\address[03]{Krasovskii Institute of Mathematics and Mechanics, S. Kovalevskaja st. 16, Yekaterinburg, 620990, Russia}
\address[04]{Novosibisk State University, Pirogova str. 2, Novosibirsk, 630090, Russia}
\address[05]{Chelyabinsk State University, Brat'ev Kashirinyh st. 129, Chelyabinsk,  454021, Russia}

\begin{abstract}
We consider the symmetric group $\mathrm{Sym}_n,\,n\geqslant 2$, generated by the set $S$ of transpositions $(1~i),\,2 \leqslant i \leqslant n$, and the Cayley graph $S_n=Cay(\mathrm{Sym}_n,S)$ called the Star graph. For any positive integers $n\geqslant 3$ and $m$ with $n > 2m$, we present a family of $PI$-eigenfunctions of $S_n$ with eigenvalue $n-m-1$. We establish a connection of these functions with the standard basis of a Specht module. In the case of largest non-principal eigenvalue $n-2$ we prove that any eigenfunction of $S_n$ can be reconstructed by its values on the second neighbourhood of a vertex.

\end{abstract}

\begin{keyword}
Cayley graph; Star graph; symmetric group; graph spectrum; eigenvalues; eigenfunctions; Jucys-Murphy elements;
\vspace{\baselineskip}
\MSC[2010] 05C25\sep 05E10\sep 05E15\sep 90B10
\end{keyword}
\end{frontmatter}

\section{Introduction}\label{sec0}

The study of eigenfunctions of graphs plays an important role in theoretical and applied research~\cite{BH12}. Eigenfunctions of graphs are related to various combinatorial structures such as perfect codes, equitable partitions, trades~\cite{KMP16,K17,V17}.

Denote by $\mathrm{Sym}_n$ the group consisting of all bijections from $\{1, 2, \ldots ,n\}$ to itself using composition $\circ$ as  multiplication. We investigate eigenfunctions of the Star graph $S_n=Cay(\mathrm{Sym}_n,S)$, $n\geqslant 2$, that is the Cayley graph on $\mathrm{Sym}_n$ generated by transpositions from the set $S=\{(1~i)~|~ 2 \leqslant i \leqslant n\}$. By the definition, the Star graph is a connected bipartite $(n-1)$--regular graph.

The spectrum of the Star graph is integral~\cite{CF12,KM75}. More precisely, for $n \geqslant 2$ and for each integer $1\leqslant k \leqslant n-1$, the values $\pm(n-k)$ are eigenvalues of $S_n$; if $n\geqslant 4$, then $0$ is an eigenvalue of $S_n$. Since the Star graph is bipartite, $\mathrm{mul}(n-k)=\mathrm{mul}(-n+k)$ for each integer $1\leqslant k \leqslant n$. Moreover, $\pm(n-1)$ are simple eigenvalues of $S_n$. By~\cite{CF12}, the integrality of the spectrum of the Star graph $S_n$ follows from studying the spectrum of the Jucys-Murphy elements $J_n$ in the algebra of the symmetric group. References on the topic see also in~\cite{R11}.

In this paper, for any positive integers $n\geqslant 3$ and $m$ with $n > 2m$, we present a family of eigenfunctions $f_{I_m}^{P_m}$ of the Star graph $S_n$ with eigenvalue $n-m-1$, where $I_m$ is a vector of $m$ pairwise different elements from the set $\{1,\ldots,n\}$ and $P_m$ is a vector of $m$ pairs of pairwise different elements from the set $\{2,\ldots,n\}$. We call these eigenfunctions as $PI$-\emph{eigenfunction} of $S_n$.

We prove that an eigenfunction of the Jucys-Murphy operator $J_n$ with eigenvalue $n-m-1$, $n > 2m$, given by a polytabloid is expressed as a sum of $PI$-eigenfunctions of $S_n$. This result follows from an observation that the adjacency matrix of $S_n$ coincides with the transformation matrix of  $J_n$ and demonstrates a remarkable connection between spectral properties of the Star graph $S_n$ and the representation theory of the symmetric group. In the case of largest non-principal eigenvalue $n-2$ we prove that any eigenfunction of $S_n$ can be reconstructed by its values on the second neighbourhood of a vertex.

The paper is organized as follows. In Section~$2$ we define a family of $PI$-eigenfunctions with eigenvalues of $S_n$ greater than $(n-2)/2$. In Section~$3$ we establish a correspondence between eigenfunctions of $J_n$ and $S_n$. In Section~$4$ we express an eigenfunction of $J_n$ with eigenvalue $n-m-1$, $n > 2m$, given by a polytabloid as a sum of $PI$-eigenfunctions of $S_n$. Finally, in Section~$5$ we prove that any eigenfunction of $S_n$ with eigenvalue $n-2$ can be reconstructed by its values on the second neighbourhood of a vertex.

\section{Family of $PI$-eigenfunctions of the Star graph $S_n$}\label{mainSec}

Let $\Gamma$ be a graph with the adjacency matrix $A$. Let $\theta$ be an eigenvalue of the matrix $A$. A function $f:V(\Gamma)\longrightarrow{\mathbb{R}}$ is called an \emph{eigenfunction} of $\Gamma$ corresponding to $\theta$ if $f\not\equiv 0$ and the equality
\begin{equation}\label{LocalCondition}
    \theta\cdot f(x)=\sum\limits_{y\in{N(x)}}f(y)
\end{equation}
holds for any its vertex $x$, where $N(x)$ is the neighborhood of $x$.

Let us define a vector $I_m=(i_1, i_2, \ldots, i_m)$ of $m$ pairwise different elements from the set $\{1,\ldots,n\}$ and a vector $P_m=((j_1,k_1), (j_2,k_2), \ldots, (j_m,k_m))$ of $2m$ pairwise different elements from the set $\{2,\ldots,n\}$ arranged into $m$ pairs. Define a function $f_{I_m}^{P_m}:\mathrm{Sym}_n \rightarrow \mathbb{R}$. For a permutation $\pi=[\pi_1\pi_2\ldots\pi_n] \in \mathrm{Sym}_n$, we put $f_{I_m}^{P_m}(\pi)=0,$ if there exists $t \in \{1,2, \ldots, m\}$ such that $\pi_{j_t} \ne i_t$ and $\pi_{k_t} \ne i_t$. If for every $t \in \{1,2, \ldots, m\}$ either $\pi_{j_t}=i_t$ or $\pi_{k_t}=i_t$, then we define a binary vector $X_\pi=(x_1, x_2, \ldots, x_m)$ as follows:
$$
x_t=
\left\{
  \begin{array}{ll}
    1, & \hbox{if $\pi_{j_t} = i_t$;} \\
    0, & \hbox{if $\pi_{k_t} = i_t$.}
  \end{array}
\right.$$
We use the vector $X_\pi$ to complete the definition of the function $f_{I_m}^{P_m}$:

\begin{equation}\label{MainFunction}
f_{I_m}^{P_m}(\pi)=
\left\{
  \begin{array}{ll}
    1, & \hbox{if $X_\pi$ contains an even number of 1s;} \\
    -1, & \hbox{if $X_\pi$ contains an odd number of 1s;}\\
    0, & \hbox{there exists $t$
such that $\pi_{j_t} \ne i_t$ and $\pi_{k_t} \ne i_t$.}
  \end{array}
\right.
\end{equation}

\begin{proposition}\label{eigenfunctions} For $n \geqslant 3$, the function $f_{I_m}^{P_m}$ is an eigenfunction with eigenvalue $n-m-1$ of the Star graph $S_n$.
\end{proposition}
\proof To show that $f_{I_m}^{P_m}$ is an eigenfunction, it is enough to verify that~$(\ref{LocalCondition})$ holds for any permutation $\pi=[\pi_1\pi_2\ldots\pi_n] \in \mathrm{Sym}_n$.

Suppose that $\pi_1=i_t$ for some $t \in \{1,2, \ldots, m\}$. Since $j_t \ne 1$ and $k_t \ne 1$, we have $f_{I_m}^{P_m}(\pi)=0$ by~(\ref{MainFunction}). Let $\sigma \in \mathrm{Sym}_n$ be a permutation such that $f_{I_m}^{P_m}(\sigma) \ne 0$ and $\sigma$ is adjacent to $\pi$. It follows from the definitions of $f_{I_m}^{P_m}$ and $S_n$ that the set of neighbours of $\pi$ with non-zero value of $f_{I_m}^{P_m}$ is either empty or equal to $\{(1~j_t)\circ\pi, (1~k_t)\circ\pi\}$. Note, that $f_{I_m}^{P_m}((1~j_t)\circ\pi)= -f_{I_m}^{P_m}((1~k_t)\circ\pi)$, which means that~(\ref{LocalCondition}) holds for the permutation $\pi$.

Suppose that $\pi_1 \ne i_t$ for any $t \in \{1, 2, \ldots, m\}$. If $f_{I_m}^{P_m}(\pi)=0$, then $\pi$ has no neighbours with non-zero value of $f_{I_m}^{P_m}$, and~(\ref{LocalCondition}) holds for the permutation $\pi$. If $f_{I_m}^{P_m}(\pi) \ne 0$, then the neighbours of $\pi$ with non-zero value of $f_{I_m}^{P_m}$ form the set $S=\{(1~s)\circ\pi~|~ s \in \{2,\ldots,n\}, \pi_s \not\in\{i_1, \ldots, i_m\}\}$ of cardinality $n-m-1$. Moreover, $f_{I_m}^{P_m}(\sigma)=f_{I_m}^{P_m}(\pi)$ for any $\sigma \in S$, which means that~(\ref{LocalCondition}) holds for the permutation $\pi$ in this case. This completes the proof. $\square$

\section{Correspondence between eigenfunctions of $J_n$ and $S_n$}\label{sec4}
In~\cite{CF12} it was pointed out that studying the spectrum of the Star graph is equivalent to studying the spectrum of the Jucys-Murphy element
\begin{equation}\label{e8}
J_n=(1~n)+(2~n)+\ldots+(n-1~n)
\end{equation}
acting on the group algebra $\mathbb{C}[\mathrm{Sym}_n]$ by left multiplication. This can be shown by the following arguments, where we adopted terminology from~\cite{Sa01}.

Denote by $S_n^{JM}$ the Cayley graph on the symmetric group $\mathrm{Sym}_n$ with the generating set $\{(1~n),(2~n),\ldots,(n-1~n)\}$, which is obviously isomorphic to the Star graph $S_n$. The Jucys-Murhpy element $J_n$ presented by~(\ref{e8}) can be considered as an adjacency operator such that for any $\pi\in\mathrm{Sym}_n$ the equality $J_n(\pi)=\sum\limits_{\sigma\in N(\pi)}\sigma$ holds, where $N(\pi)$ is the neighbourhood of $\pi$ in $S_n^{JM}$. For any element $v=\sum\limits_{\pi\in\mathrm{Sym}_n}y_\pi\pi$ from $\mathbb{C}[\mathrm{Sym}_n]$, we define a function $f_v:\mathrm{Sym}_n\rightarrow \mathbb{C}$ such that for any $\pi\in\mathrm{Sym}_n$ the equality $f_v(\pi)=y_\pi$ holds. The mapping $\psi:v\rightarrow f_v$ gives a bijection between the elements of $\mathbb{C}[\mathrm{Sym}_n]$ and the complex-valued functions on $\mathrm{Sym}_n$. The following lemma shows a correspondence between eigenfunctions of $J_n$ and $S_n^{JM}$.
\begin{lemma}\label{EigenJucysandEigenStar}
Let $v=\sum\limits_{\pi\in \mathrm{Sym}_n}y_\pi\pi \in \mathbb{C}[\mathrm{Sym}_n]$ be an eigenfunction of $J_n$ with eigenvalue $\theta$ and let $y_\pi$'s be real. Then $\psi(v)=f_v$ is an eigenfunction of $S_n^{JM}$ with eigenvalue $\theta$.
\end{lemma}
\proof It follows from the fact that the adjacency matrix of $S_n^{JM}$ and the transformation matrix of the operator $J_n$
acting on $\mathbb{C}[\mathrm{Sym}_n]$ coincide. $\square$

\medskip
Let $\lambda$ be a shape with $n$ cells. For a tableau $t$ of shape $\lambda$, the $\lambda$-\emph{tabloid} $\{t\}$ is the set of all tableaux of shape $\lambda$ that can be obtained from $t$ by permutations of elements in rows. Let $M^\lambda=\mathbb{C}\{\{t_1\},\ldots,\{t_k\}\}$ be the permutation module corresponding to $\lambda$, where $\{t_1\},\ldots,\{t_k\}$ is a complete list of $\lambda$-tabloids. We establish a correspondence between eigenfunctions of $J_n$ acting on the permutation module $M^\lambda$ corresponding to shape $\lambda$ and eigenfunctions of $J_n$ acting on the group algebra $\mathbb{C}[\mathrm{Sym}_n]$.

Let $t$ be a tableau of shape $\lambda$ and $C_t$ be the column-stabilizer of $t$.
\begin{lemma}[\cite{Sa01}, Lemma 2.3.3(2)]\label{phicolstab}
For any $\pi\in\mathrm{Sym}_n$, the equality $C_{\pi(t)}=\pi^{-1}\circ C_t\circ \pi$ holds.
\end{lemma}

For any tableau $t$ of shape $\lambda$, the element $\mathbf{e}_t=\sum\limits_{\sigma\in C_t}sgn(\sigma)\{\sigma(t)\}$ is called
the \emph{polytabloid}. The submodule $S^\lambda$ of $M^\lambda$ spanned by all polytabloids $\mathbf{e}_t$, where $t$ is a tableau of shape $\lambda$, is called the \emph{Specht module} associated with $\lambda$. It is well known that the set of \emph{standard polytabloids} $$\{\mathbf{e}_t:\text{$t$ is a standard tableau of shape $\lambda$}\}$$ is a basis for $S^\lambda$. Moreover, the Specht modules form a complete list of irreducible $\mathrm{Sym}_n$-modules over the complex field. On the other hand, there exists a basis for $S^\lambda$ consisting of eigenfunctions of $J_n$ (see~\cite[Theorem 1.1]{CF12} and~\cite{J74}). Since the set of standard polytabloids is a basis for $S^\lambda$, each eigenfunction of $J_n$ in $S^\lambda$ is a linear combination of standard polytabloids.

Let $id_\lambda$ be the standard tableau of shape $\lambda$ whose rows consist of the consecutive elements. Let $T_\lambda$ be the set of all tableaux of shape $\lambda$. For any tableau $t\in T_\lambda$, denote by $\tau_t$ the permutation defined be the equation
\begin{equation}\label{tau}
\tau_t(t)=id_\lambda,
\end{equation}
where $\tau_t$ acts on $t$ by replacing the values of the cells of $t$ by their images under the permutation mapping $\tau_t$.

\begin{lemma}\label{tausigma}
For any $t\in T_\lambda$ and $\sigma\in\mathrm{Sym}_n$, the equality $\tau_{\sigma(t)}=\sigma^{-1}\circ\tau_t$ holds.
\end{lemma}
\proof
$(\sigma^{-1}\circ\tau_t)(\sigma(t))=\tau_t(\sigma^{-1}(\sigma(t)))=id_\lambda. \square$

\medskip
Let us define a linear mapping $\phi:M^\lambda\rightarrow \mathbb{C}[\mathrm{Sym}_n]$. Since the set of all $\lambda$-tabloids is
a basis for $M_\lambda$, it is enough to define images for $\lambda$-tabloids. For any $\lambda$-tabloid $\{t\}$, where $t \in T_\lambda$, we put $\phi(\{t\})=\sum\limits_{t'\in\{t\}}\tau_{t'}$. Since the $\lambda$-tabloids $\{t_1\},\ldots,\{t_k\}$ form a partition of $T_\lambda$, the mapping $\phi$ is an isomorphic embedding of $M^\lambda$ into $\mathbb{C}[\mathrm{Sym}_n]$.

\begin{lemma}\label{phisigma}
For any $\lambda$-tabloid $\{t\}$ and $\sigma\in\mathrm{Sym}_n$, the equality $\phi(\{\sigma(t)\})=\sigma^{-1}\circ\phi(\{t\})$ holds.
\end{lemma}
\proof
$$\phi(\{\sigma(t)\})=\sum\limits_{t'\in\{\sigma(t)\}}\tau_{t'}=\sum\limits_{t'\in\{t\}}\tau_{\sigma(t')}=$$
(by Lemma \ref{tausigma})
$$=\sigma^{-1}\circ\sum\limits_{t'\in\{t\}}\tau_{t'}=\sigma^{-1}\circ\phi(\{t\}).$$ $\square$

\medskip
Lemma \ref{phipolytabloid} immediately follows from Lemma \ref{phisigma}.
\begin{lemma}\label{phipolytabloid}
For any polytabloid $\mathbf{e}_t$, the equality $\phi(\mathbf{e}_t)=\sum\limits_{\pi\in C_t}sgn(\pi)\pi^{-1}\circ
\phi(\{t\})
$ holds.
\end{lemma}

Let us put $\sigma_i=(i~n)$ up to the end of this section. Then $J_n=\sum\limits_{i=1}^{n-1}\sigma_i$ holds.

\begin{lemma}\label{phiandJcommute}
For any polytabloid $\mathbf{e}_t$, the equality $\phi(J_n(\mathbf{e}_t))=J_n(\phi(\mathbf{e}_t))$ holds.
\end{lemma}
\proof
$$\phi(J_n(\mathbf{e}_t))=\phi(\sum\limits_{i=1}^{n-1}\sigma_i(\mathbf{e}_t))=\phi(\sum\limits_{i=1}^{n-1}\mathbf{e}_{\sigma_i(t)})=$$
(since $\phi$ is linear)
$$=\sum\limits_{i=1}^{n-1}\phi(\mathbf{e}_{\sigma_i(t)})=$$
(by Lemmas \ref{phipolytabloid} and \ref{phisigma})
$$=\sum\limits_{i=1}^{n-1}\sum\limits_{\pi\in C_{\sigma_i(t)}}sgn(\pi)\pi^{-1}\circ\sigma_i^{-1}\circ \phi(\{t\})=$$
(by Lemma \ref{phicolstab})
$$=\sum\limits_{i=1}^{n-1}\sigma_i^{-1}\circ\sum\limits_{\pi\in C_t}sgn(\pi)\phi(\{\pi(t)\})=\sum\limits_{i=1}^{n-1}\sigma_i^{-1}\circ\phi(\mathbf{e}_t)=J_n(\phi(\mathbf{e}_t)).$$
$\square$

\begin{proposition}\label{equiv}
Let $\mathbf{v}\in S^\lambda$ be an eigenfunction of the operator $J_n:M^\lambda\rightarrow M^\lambda$ with eigenvalue $\theta$. Then $\phi(\mathbf{v})$ is an eigenfunction of the operator $J_n:\mathbb{C}[\mathrm{Sym}_n]\rightarrow\mathbb{C}[\mathrm{Sym}_n]$ with  eigenvalue $\theta$.
\end{proposition}
\proof Let $\mathbf{v}=\sum\limits_{i}a_i\mathbf{e}_{t_i}$, then we have
$$J_n(\phi(\mathbf{v}))=J_n(\phi(\sum\limits_{i}a_i\mathbf{e}_{t_i}))=\sum\limits_{i}a_iJ_n(\phi(\mathbf{e}_{t_i}))=$$
(by Lemma~\ref{phiandJcommute})
$$=\sum\limits_{i}a_i\phi(J_n(\mathbf{e}_{t_i}))=\phi(J_n(\sum\limits_{i}a_i\mathbf{e}_{t_i}))=\phi(J_n(\mathbf{v}))=
\phi(\theta \mathbf{v})=\theta \phi(\mathbf{v}).$$ $\square$

\medskip
Proposition~\ref{equiv} and Lemma~\ref{EigenJucysandEigenStar} give a way to obtain eigenfunctions of $S_n^{JM}$ from eigenfunctions of $J_n$, where $J_n$ acts on $M^\lambda$. Due to the isomorphism between $S_n$ and $S_n^{JM}$, if we put
$$\{j_1, k_1, j_2, k_2, \ldots, j_m, k_m\} \subset \{1, 2, 3, \ldots, n-1\}$$ in setting~(\ref{MainFunction}), then a $PI$-eigenfunction $f_{I_m}^{P_m}$ of $S_n$ becomes a $PI$-eigenfunction of the graph $S_n^{JM}$.

The following example illustrates a certain connection between $PI$-eigenfunctions of $S_n^{JM}$ and eigenfunctions of $J_n$ given by standard polytabloids.

\textbf{Example 1.} Let us take the partition $\lambda=(n-1,1)$. For any $i \in \{2,\ldots,n-1\}$, consider the standard tableau
{
\ytableausetup{centertableaux}
$
t_i=
\begin{ytableau}
1 & \ldots & n \\
i
\end{ytableau}
$
}, its corresponding polytabloid $\mathbf{e}_{t_i}$ and put
{
\ytableausetup{centertableaux}
$
t_1=
\begin{ytableau}
2 & \ldots & n \\
1
\end{ytableau}
$
},
{
\ytableausetup{centertableaux}
$
t_n=
\begin{ytableau}
1 & \ldots & \scriptstyle n-1 \\
n
\end{ytableau}
$
}.

Then we have the following equalities:
$$\mathbf{e}_{t_i} = \{t_i\} - \{t_1\},$$
$$\sigma_\ell(\mathbf{e}_{t_i}) = \mathbf{e}_{t_i}, \text{for any $\ell \in \{1,\ldots,n-1\}\setminus\{1,i\}$ },$$
$$\sigma_1(\mathbf{e}_{t_i}) = \{t_i\} - \{t_n\},$$
$$\sigma_i(\mathbf{e}_{t_i}) = \{t_n\} - \{t_1\},$$
$$J_n(\mathbf{e}_{t_i}) = \sigma_1(\mathbf{e}_{t_i}) + \sigma_i(\mathbf{e}_{t_i}) + \sum\limits_{\ell = 2,~\ell \ne i}^{n-1}\sigma_\ell(\mathbf{e}_{t_i}) =
(n-2)\mathbf{e}_{t_i}.
$$
By Proposition \ref{equiv}, the elements $\phi(\mathbf{e}_{t_2}), \ldots, \phi(\mathbf{e}_{t_{n-1}})$ belonging to $\mathbb{C}[\mathrm{Sym}_n]$ are eigenfunctions of $J_n$ corresponding to the eigenvalue $n-2$. Note, that for any $i \in \{2,\ldots,n-1\}$ the equality $\phi(\mathbf{e}_{t_i})=\sum\limits_{t'_i\in \{t_i\} }\tau_{t'_i} - \sum\limits_{t'_1\in \{t_1\}}\tau_{t'_1}$ holds.
Then we have
\begin{equation}\label{EqSimplyMainFunction}
f_{\phi(\mathbf{e}_{t_i})}(\pi)=
\left\{
  \begin{array}{ll}
    1, & \hbox{if $\pi_i = n$;} \\
    -1, & \hbox{if $\pi_1 = n$;}\\
     0, & \hbox{otherwise,}
  \end{array}
\right.
\end{equation}
where $\pi \in \mathrm{Sym}_n$. $\square$

\medskip
We have seen in Example 1 that the function (\ref{EqSimplyMainFunction}) coincides with the $PI$-eigenfunction $f^{i,1}_n$ of $S_n^{JM}$.

Now we investigate the polytabloids $\mathbf{e}_t$ corresponding to the tableaux $t$ defined below. For any integer $n \geqslant 3$ denote by $\mathcal{P}(n)$ the set of partitions of $n$. Let $\lambda \in \mathcal{P}(n)$ be a partition $(\lambda_1, \lambda_2, \ldots, \lambda_s)$, where $s \geqslant 2$, $\lambda_1 > \lambda_2$ and $\lambda_i \geqslant \lambda_{i+1}$ for any $i \in \{2, \ldots, s-1\}$. Put $m = \lambda_2 + \ldots + \lambda_s$. For shape $\lambda$, denote by $\mu = (\mu_1, \mu_2, \ldots, \mu_{\lambda_1})$ the conjugate of $\lambda$, where $\mu_i$ is the length of $i$th column of $\lambda$, and let $k$ be the number of columns of length greater than one. In this setting $m$ is the number of cells in all rows of $\lambda$ but the first.

Let $t$ be a standard tableau of shape $\lambda$ with $n$ placed at its upper right cell. For the tableau $t$, denote by $X_t$ the set of $m + k$ elements in $k$ columns of length greater than one and by $Y_t$ the set of $n-m-k-1$ elements in all columns of length one but the last.

In Proposition \ref{PolytabloidIsAnEigenfunction}, we prove that the polytabloid $\mathbf{e}_t$ is an eigenfunction of the Jucys-Murphy operator $J_n$. The similar proof can be found in~\cite[Theorem 3.7]{FOW85}.
\begin{proposition}\label{PolytabloidIsAnEigenfunction}
The polytabloid $\mathbf{e}_t$ is an eigenfunction of $J_n$ with eigenvalue $n-m-1$.
\end{proposition}
\proof
Put $J_n^{X_t}=\sum\limits_{i\in X_t}\sigma_i$ and $J_n^{Y_t}=\sum\limits_{\ell\in Y_t}\sigma_\ell$, then we have
$J_n(\mathbf{e}_t)=J_n^{X_t}(\mathbf{e}_t)+J_n^{Y_t}(\mathbf{e}_t)=\sum\limits_{i \in X_t}\sigma_i(\mathbf{e}_t)+\sum\limits_{\ell \in Y_t}\sigma_l(\mathbf{e}_t)$. Since for any $\ell \in Y_t$ the equality $\sigma_\ell(\mathbf{e}_t)=\mathbf{e}_t$ holds, we have $$J_n^{Y_t}(\mathbf{e}_t) = |Y_t|\mathbf{e}_t = (n-m-k-1)\mathbf{e}_t.$$

Now, it is enough to show that $J_n^{X_t}(\mathbf{e}_t) = k\mathbf{e}_t$ holds. By definition of the polytabloid, $J_n^{X_t}(\mathbf{e}_t)$ is equal to the sum $\sum\limits_{i \in X_t}\sum\limits_{\pi \in C_t}sgn(\pi)\{\sigma_i(\pi(t))\}$, where each of the $|X_t||C_t|$ summands in this sum is uniquely determined by the pair $i$ and $\pi$.

We say that an element from \{1, \ldots, n\} belongs to a row of a $\lambda$-tabloid, if this element belongs to the corresponding row of each tableau of this $\lambda$-tabloid.

Firstly, we prove that the $\lambda$-tabloids with no $n$ in the first row are annihilated in $J_n^{X_t}(\mathbf{e}_t)$. Let $\{\sigma_{i_2}({\pi(t)}\})$ be an arbitrary $\lambda$-tabloid in $J_n^{X_t}(\mathbf{e}_t)$ with no $n$ in the first row for some $i_2$ from $X_t$. Then $i_2$ is an element in a non-first row in the tableau $\pi(t)$. Denote by $i_1$ the element of $\pi(t)$ belonging to the first row and the same column as $i_2$. Consider the permutation $\pi\circ (i_1~i_2) \in C_t$ and denote it by $\pi'$. Then the tableaux $\pi(t)$ and $\pi'(t)$ are only differed by the cells containing either $i_1$ or $i_2$. We have the same property for the tableaux $\sigma_{i_2}(\pi(t))$ and $\sigma_{i_1}(\pi'(t))$. Moreover, the elements $i_1,i_2$ belong to the first row of each of the tableaux $\sigma_{i_2}(\pi(t))$ and $\sigma_{i_1}(\pi'(t))$. Hence, the $\lambda$-tabloids $\{\sigma_{i_2}(\pi(t))\}$ and $\{\sigma_{i_1}(\pi'(t))\}$ coincide. In the sum $J_n^{X_t}(\mathbf{e}_t)$, the $\lambda$-tabloids $\{\sigma_{i_2}(\pi(t))\}$ and $\{\sigma_{i_1}(\pi'(t))\}$ are taken with different signs since the permutations $\pi$ and $\pi'$ have different parities.

Now we show that each $\lambda$-tabloid with $n$ in the first row appears in $J_n^{X_t}(\mathbf{e}_t)$ exactly $k$ times. Since the tableau $t$ has exactly $k$ elements from $X_t$ in the first row and $C_t$ permutes only elements of $X_t$, then for any $\pi \in C_t$ the tableau $\pi(t)$ and, as a result, the $\lambda$-tabloid $\{\pi(t)\}$ have exactly $k$ elements from $X_t$ in the first row. Denote by $i_1, \ldots, i_k$ these $k$ elements of $X_t$ in the first row of $\{\pi(t)\}$. The only transpositions of $J_n^{X_t}$ that permute elements in the first row of $\{\pi(t)\}$ are $\sigma_{i_1}, \ldots, \sigma_{i_k}$. Hence, the transpositions of $J_n^{X_t}$ that stabilize the $\lambda$-tabloid $\{\pi(t)\}$ are the same. Thus, each of the $\lambda$-tabloids $\{\pi(t)\}$ forming the polytabloid $\mathbf{e}_t$ appears in $J_n^{X_t}(\mathbf{e}_t)$ with the sign of $\pi$ exactly $k$ times, and $J_n^{X_t}(\mathbf{e}_t) = k\mathbf{e}_t$ holds. The proposition is proved. $\square$

\textbf{Example 2.} We illustrate Proposition \ref{PolytabloidIsAnEigenfunction} for the partition $\lambda=(n-2,2)$. For any pair $i,j \in \{1,\ldots,n\}, i \ne j$,
put
{
\ytableausetup{centertableaux}
$
t_{ij}=
\begin{ytableau}
\scriptstyle i_1 & \scriptstyle i_2 &\ldots & \scriptstyle i_{n-2} \\
i & j
\end{ytableau}
$
}, where $\{i_1, \ldots, i_{n-2}\} = \{1,2,\ldots,n\} \setminus \{i,j\}$ and $i_1 < i_2 < \ldots < i_n$ hold.

Take any $i,j \in \{2,\dots, n-1\}, i < j,$ and consider the standard polytabloid $\mathbf{e}_{t_{ij}}$
corresponding to the tableau
{
\ytableausetup{centertableaux}
$
t_{ij}=
\begin{ytableau}
1 & i_2 & \ldots & n \\
i & j
\end{ytableau}
$
}.
In this setting we have $X_t = \{1,i_2,i,j\}$, $Y_t = \{1,\ldots,n-1\}\setminus X_t$, and the following equalities hold:
$$\mathbf{e}_{t_{ij}} = \{t_{ij}\} + \{t_{1i_2}\} - \{t_{1j}\} - \{t_{ii_2}\},$$
$$\sigma_\ell(\mathbf{e}_{t_{ij}}) = \mathbf{e}_{t_{ij}}, \text{for any $\ell \in \{1,\ldots,n-1\}\setminus\{1,i_2, i, j\}$},$$
$$\sigma_1(\mathbf{e}_{t_{ij}}) = \{t_{ij}\} + \{t_{ni_2}\}- \{t_{nj}\} - \{t_{ii_2}\},$$
$$\sigma_{i_2}(\mathbf{e}_{t_{ij}}) = \{t_{ij}\} + \{t_{1n}\}- \{t_{1j}\} - \{t_{in}\},$$
$$\sigma_i(\mathbf{e}_{t_{ij}}) = \{t_{nj}\} + \{t_{1i_2}\} - \{t_{1j}\} - \{t_{ni_2}\},$$
$$\sigma_j(\mathbf{e}_{t_{ij}}) = \{t_{in}\} + \{t_{1i_2}\}- \{t_{1n}\} - \{t_{ii_2}\},$$
$$J_n(\mathbf{e}_{t_{ij}}) = J_n^{X_t}(\mathbf{e}_{t_{ij}}) + J_n^{Y_t}(\mathbf{e}_{t_{ij}})=
2\mathbf{e}_{t_{ij}} + (n-5)\mathbf{e}_{t_{ij}}
= (n-3)\mathbf{e}_{t_{ij}},
$$
where
$$J_n^{X_t}(\mathbf{e}_{t_{ij}}) = \sigma_1(\mathbf{e}_{t_{ij}}) + \sigma_{i_2}(\mathbf{e}_{t_{ij}}) + \sigma_i(\mathbf{e}_{t_{ij}}) + \sigma_j(\mathbf{e}_{t_{ij}}),$$
$$
J_n^{Y_t}(\mathbf{e}_{t_{ij}})= \sum\limits_{\ell = 2,~\ell \ne i_2,i,j}^{n-1}\sigma_\ell(\mathbf{e}_{t_{ij}}).
$$

\medskip
Note that the equality
\begin{equation}\label{PolyIsASumOfTwo}
f_{\phi(\mathbf{e}_{t_{ij}})} = f^{(i,1),(j,i_2)}_{n-1,n} + f^{(j,i_2),(i,1)}_{n-1,n}.
\end{equation} holds.
$\square$

In the next section we generalize Equality (\ref{PolyIsASumOfTwo}).

\section{Decomposition of eigenfunctions into $PI$-eigenfunctions }\label{PolytabloidFunctionInTermsOfOurFunctions}

In this section, for a standard tableau $t$ of shape $\lambda$ and the polytabloid $\mathbf{e}_t$ from Proposition~\ref{PolytabloidIsAnEigenfunction}, we express the eigenfunction  $f_{ \phi( \mathbf{e}_{t} ) }$ as a sum of certain $PI$-eigenfunctions of $S_n^{JM}$ in the case when $n > 2m$.

Put $X_\lambda=\{1, 2, \ldots, n-m\}$. The subgroup $Sym(X_\lambda)$ of $Sym_n$ is the pointwise stabilizer of the set $\{n-m+1, n-m+2, \ldots, n\}$.  For any left coset of $Sym(X_\lambda)$ in $Sym_n$, there exists uniquely determined $m$-element sequence $z_1,z_2, \ldots, z_m$ such that this coset consists of all permutations of the form
$$
\left(
  \begin{array}{cccccccc}
    * & * & \ldots & * & z_1 & z_2 & \ldots & z_m \\
    1 & 2 & \ldots & n-m & n-m+1 & n-m+2 & \ldots & n \\
  \end{array}
\right).
$$
The sequence $(z_1,z_2, \ldots, z_m)$ is called the \emph{code} of the coset, and the coset is called $(z_1,z_2, \ldots, z_m)$-\emph{coset} and is denoted by $\overline{(z_1,z_2, \ldots, z_m)}$. The following statement is straightforward.
\begin{lemma}\label{CosetsAreEqual}
Two cosets coincide if and only if their codes are equal.
\end{lemma}

For a tableau $t$ of shape $\lambda$, denote by $[t]$ the set of all tableaux of shape $\lambda$ that have the same $2$nd, $\ldots$, $s$th rows as $t$ with $|[t]| = (n-m)!$. Put $\tau_{[t]} = \{\tau_{t'}~|~ t' \in [t]\}$, where $\tau_{t'}$ is given by (\ref{tau}).

\begin{lemma}\label{Cosets}
Let $t$ be a tableau. The set $\tau_{[t]}$ is the $(z_1,z_2, \ldots, z_m)$-coset if and only if $z_1,z_2, \ldots, z_m$ is the sequence obtained by concatenation of $2$nd, $\ldots$, $s$th rows of $t$.
\end{lemma}
\proof It follows directly from (\ref{tau}) and the fact that the rows of the tableau $id_\lambda$ consist of consecutive elements. $\square$

\medskip
The set of tableaux $[t]$ is called the $(z_1,z_2, \ldots, z_m)$-\emph{coset} (or just a \emph{coset}) if $\tau_{[t]}$
is the $(z_1,z_2, \ldots, z_m)$-\emph{coset}.

\begin{lemma}\label{ActionOnCosets}
Let $[t]$ be a $(z_1,z_2, \ldots, z_m)$-coset and $\pi$ be a permutation. Then equality $[\pi(t)] = \pi([t])$ holds.
\end{lemma}
\proof It follows from the fact that $\pi$ acts cellwise on a tableau. $\square$

\medskip
For any $i \in \{1, \ldots, s\}$ and $j \in \{1, \ldots, k\}$, denote by $R_t(i)$ and $C_t(j)$ the symmetric groups on the elements of $i$th row and $j$th column of the tableau $t$, respectively. Then we have
$$R_t = R_t(1) \times R_t(2) \times \ldots \times R_t(s),$$
$$C_t = C_t(1) \times C_t(2) \times \ldots \times C_t(k),$$
where $R_t$ and $C_t$ are the row-stabilizer and the column-stabilizer of $t$, respectively.

\begin{lemma}\label{TabloidIntoCosets}
Let $t$ be a tableau of shape $\lambda$. 
The $\lambda$-tabloid $\{t\}$ can be partitioned into cosets as follows:
$$
\{t\} = \bigcup\limits_{\sigma \in R_t(2) \times \ldots \times R_t(s)} [\sigma(t)].
$$
\end{lemma}
\proof It follows immediately from the definition of a $\lambda$-tabloid. $\square$

\medskip
For a function $f:Sym_n\rightarrow \mathbb{R}$, put
$$Supp^+_{T_\lambda}(f) = \{t' \in T_\lambda ~|~f(\tau_{t'}) > 0\},$$
$$Supp^-_{T_\lambda}(f) = \{t' \in T_\lambda ~|~f(\tau_{t'}) < 0\},$$
where $T_\lambda$ is the set of all tableaux of shape $\lambda$. The set $Supp^+_{T_\lambda}(f) \cup Supp^-_{T_\lambda}(f)$
is called the \emph{support} of the function $f$. A function $f:Sym_n\rightarrow \mathbb{R}$ is called a $(0,-1,1)$-\emph{function}
if all values of $f$ are taken from the set $\{0,-1,1\}$.

\begin{lemma}\label{SuppPolytabloid}
The function $f_{\phi(\mathbf{e}_t)}$ is a $(0,-1,1)$-function, where:
$$Supp^+_{T_\lambda}(f_{\phi(\mathbf{e}_t)}) = \bigcup\limits_{\substack{\pi \in C_t,\\\pi~\text{is even}} }
\bigcup\limits_{\sigma \in R_{\pi(t)}(2) \times \ldots \times R_{\pi(t)}(s)} [\sigma(\pi(t))]
;$$
$$Supp^-_{T_\lambda}(f_{\phi(\mathbf{e}_t)}) = \bigcup\limits_{\substack{\pi \in C_t,\\\pi~\text{is odd}} }
\bigcup\limits_{\sigma \in R_{\pi(t)}(2) \times \ldots \times R_{\pi(t)}(s)} [\sigma(\pi(t))]
.$$
\end{lemma}
\proof
It follows from the definition of a polytabloid and Lemma \ref{TabloidIntoCosets}.
$\square$

\medskip
By Lemma~\ref{SuppPolytabloid}, the support of $f_{\phi(\mathbf{e}_t)}$ is a disjoint union of all cosets, whose code contains the elements from $X_t$ only, where $X_t$ is the set of $m + k$ elements in $k$ columns of length greater than one in the tableau $t$.

Remind that our main goal is to express the eigenfunction  $f_{\phi(\mathbf{e}_{t}) }$ as a sum of certain $PI$-eigenfunctions
of $S_n^{JM}$. The key argument of our proof is the fact that the supports of  $f_{\phi(\mathbf{e}_{t})}$ and a $PI$-eigenfunction can be partitioned into cosets of the subgroup $Sym(X_\lambda)$. By Lemma~\ref{SuppPolytabloid}, we have already obtained this for the eigenfunction $f_{\phi(\mathbf{e}_{t})}$.

In fact, the support of a $PI$-eigenfunction $f_{I_m}^{P_m}$ can also be partitioned into cosets of a pointwise stabilizer of $I_m$. We are especially interested in $PI$-eigenfunctions $f_{I_m}^{P_m}$ of $S_n^{JM}$ defined by vectors $I_m=(n-m+1,n-m+2,\ldots,n)$ and $P_m$, where $P_m$ is specified by introducing additional notation below. For them, we build such partition in Lemma~\ref{SupPICoset}. We assume that the tableau $t$ is represented in the following form:
$$
t
=
\ytableausetup{centertableaux, boxsize=1.8em}
\begin{ytableau}
x_{11} & x_{21} & \ldots  & \ldots  & \ldots  & \ldots  & \ldots  & \ldots & x_{k1} &  \leftarrow & Y_t & \rightarrow & n    \\
x_{12} & x_{22}  & \ldots & \ldots & \ldots  & \ldots & \ldots  & \ldots & x_{k2} \\
x_{13} & x_{23}  & \ldots & \ldots & \ldots  & \ldots & x_{\lambda_33} \\
\ldots & \ldots & \ldots & \ldots & \ldots \\
x_{1s} & x_{2s} & \ldots & x_{\lambda_ss} \\
\end{ytableau},
$$
where $x_{ji}$'s form the set $X_t$ for all appropriate $j \in \{1, \ldots, k\}$,  $i \in \{1, \ldots, s\}$, and $Y_t$ is the set of $n-m-k-1$ elements in all columns of length one but the last. Note that the picture above does not represent all possibilities. Since $n > 2m$, we have
\begin{equation}\label{Ineq}
|Y_t| \geqslant m-k.
\end{equation}
Put
$$CA_t = CA_t(1) \times CA_t(2) \times \ldots \times CA_t(k),$$
where $CA_t(j)$ denotes the subgroup of even permutations in $C_t(j)$. For any permutation $\pi \in CA_t$ there exists uniquely determined permutations $\pi_1 \in CA_t(1), \ldots, \pi_k \in CA_t(k)$ such that the equality $\pi = \pi_1 \circ \ldots \circ \pi_k$ holds, where
$$
\pi_j =
\left(
  \begin{array}{cccc}
    x_{j1} & x_{j2} & \ldots & x_{j\mu_j} \\
    \pi_{j1} & \pi_{j2} & \ldots & \pi_{j\mu_j} \\
  \end{array}
\right).
$$
Then the tableau $\pi(t)$ is presented as follows:

\begin{equation}\label{Pit}
\pi(t)
=
\ytableausetup{centertableaux, boxsize=1.8em}
\begin{ytableau}
\pi_{11} & \pi_{21} & \ldots  & \ldots  & \ldots & \ldots & \ldots & \ldots & \pi_{k1} & \leftarrow & Y_t & \rightarrow & n    \\
\pi_{12} & \pi_{22}  & \ldots & \ldots & \ldots & \ldots & \ldots & \ldots & \pi_{k2} \\
\pi_{13} & \pi_{23}  & \ldots & \ldots & \ldots & \ldots & \pi_{\lambda_33} \\
\ldots & \ldots & \ldots & \ldots & \ldots \\
\pi_{1s} & \pi_{2s} & \ldots & \pi_{\lambda_ss} \\
\end{ytableau}
\end{equation}

\medskip
Now we introduce the following two sequences of length $k$ and $m-k$. The first sequence consists of $k$ pairs of elements belonging to the columns of the rectangle of size $2\times k$ in the top-left corner of the tableau $\pi(t)$ given by (\ref{Pit}), and is presented by
\begin{equation}\label{PPi1}
(\pi_{11},\pi_{12}), (\pi_{21},\pi_{22}) ,\ldots,(\pi_{k1},\pi_{k2}).
\end{equation}
The second sequence consists of $m-k$ pairs as follows:
\begin{equation}\label{PPi2}
(y_{13},\pi_{13}), (y_{23},\pi_{23}), \ldots, (y_{\lambda_ss},\pi_{\lambda_ss}),
\end{equation}
where the pairs are obtained by matching elements from the sequences below:
$$y_{13}, y_{23}, \ldots, y_{\lambda_33}, y_{14}, \ldots, y_{\lambda_44}, \ldots
y_{1s}, \ldots, y_{\lambda_ss},$$
$$\pi_{13}, \pi_{23}, \ldots, \pi_{\lambda_33}, \pi_{14}, \ldots, \pi_{\lambda_44}, \ldots
\pi_{1s}, \ldots, \pi_{\lambda_ss},$$
where $y_{13},\ldots,y_{\lambda_ss}$ are pairwise distinct $m-k$ elements from $Y_t$ (we can always do this by~(\ref{Ineq})). Take any $\sigma \in R_{t}(2) \times \ldots \times R_{t}(s)$ and, for any $\pi \in CA_{\sigma(t)}$, we define the vector $P_\pi$ of length $m$ by concatenation (\ref{PPi1}) and (\ref{PPi2}) as follows:
$$P_\pi = ( (\pi_{11},\pi_{12}), (\pi_{21},\pi_{22}) ,\ldots,(\pi_{k1},\pi_{k2}),
(y_{13},\pi_{13}), (y_{23},\pi_{23}), \ldots, (y_{\lambda_ss},\pi_{\lambda_ss})).$$

Now, for the $PI$-eigenfunction $f^{P_\pi}_{(n-m+1,\ldots,n)}$, we give explicitly the partition of its support into cosets. For a pair $(a,b)$, where $a,b \in \{1,\ldots,n\}$, and any $\delta \in \{0,1\}$, put
$$
(a,b)^\delta =
\left\{
  \begin{array}{ll}
    a, & \hbox{if $\delta = 1$;} \\
    b, & \hbox{if $\delta = 0$.}
  \end{array}
\right.
$$
For any binary vector
$\Omega = ( \varepsilon_1, \varepsilon_2 ,\ldots, \varepsilon_k, \omega_{13}, \omega_{23}, \ldots, \omega_{\lambda_ss})$
of length $m$ and any $P_{\pi}$, put
$$P^\Omega_{\pi} = ( (\pi_{11},\pi_{12})^{\varepsilon_1},\ldots,(\pi_{k1},\pi_{k2})^{\varepsilon_k},
(y_{13},\pi_{13})^{\omega_{13}},  \ldots, (y_{\lambda_ss},\pi_{\lambda_ss})^{\omega_{\lambda_ss}}),$$
which has the following property.
\begin{lemma}\label{XandYCosets}
There exists $2^k$ values of $\Omega$ such that $P^\Omega_{\pi}$ has no elements from $Y_t$ and $2^m-2^k$ values  of $\Omega$ such that $P^\Omega_{\pi}$ has at least one element from $Y_t$.
\end{lemma}
\proof The vector $P^\Omega_{\pi}$ has no elements from $Y_t$ if and only if the vector $\Omega$ has the form
$(\underbrace{\varepsilon_1,\ldots,\varepsilon_k}_k,\underbrace{0,\ldots,0}_{m-k})$ for some $\varepsilon_1,\ldots,\varepsilon_k$. $\square$

The vector $P_\pi^\Omega$ considered as a code represents the coset $\overline{P^\Omega_{\pi}}$ which consists of the following tableaux of shape $\lambda$:
$$
\begin{array}{|@{}l@{}|@{}l@{}|@{}c@{}|@{}c@{}|@{}c@{}|@{}c@{}|@{}c@{}|@{}c@{}|@{}c@{}|@{}c@{}|@{}c@{}|}
\hline
   \multicolumn{1}{|c|}{\scriptstyle \ast}& \multicolumn{1}{c|}{\scriptstyle \ast} & \scriptstyle \dots  & \scriptstyle \ast & \scriptstyle \dots  & \scriptstyle \dots & \scriptstyle \ast &\scriptstyle \dots & \scriptstyle \ast & \scriptstyle \dots  & \scriptstyle \ast  \\
  \hline
 \scriptstyle (\pi_{11},\pi_{12})^{\varepsilon_1} & \scriptstyle  (\pi_{21},\pi_{22})^{\varepsilon_2} & \scriptstyle \dots  & \scriptstyle \dots & \scriptstyle \dots  & \scriptstyle \dots & \scriptstyle \dots & \scriptstyle \dots & \scriptstyle (\pi_{k1},\pi_{k2})^{\varepsilon_k} \\
 \cline{1-9}
\scriptstyle (y_{13},\pi_{13})^{\omega_{13}} & \scriptstyle (y_{23},\pi_{23})^{\omega_{23}} & \scriptstyle \dots  & \scriptstyle \dots  & \scriptstyle \dots & \scriptstyle \dots  & \scriptstyle (y_{\lambda_33},\pi_{\lambda_33})^{\omega_{\lambda_33}} \\
\cline{1-7}
\multicolumn{1}{|c|}{\scriptstyle \dots} & \multicolumn{1}{c|}{\scriptstyle \dots} &\scriptstyle \dots & \scriptstyle \dots & \scriptstyle \dots\\
\cline{1-5}
\scriptstyle (y_{1s},\pi_{1s})^{\omega_{1s}} & \scriptstyle (y_{2s},\pi_{2s})^{\omega_{2s}} & \scriptstyle \dots   & \scriptstyle (y_{\lambda_ss},\pi_{\lambda_ss})^{\omega_{\lambda_ss}} \\
\cline{1-4}
\end{array}
$$

\begin{lemma}\label{SupPICoset} For any  $\sigma \in R_{t}(2) \times \ldots \times R_{t}(s)$ and $\pi \in CA_{\sigma(t)}$,
the following equalities hold:\\
{\rm(1)}~$Supp^+_{T_\lambda}(f^{P_\pi}_{(n-m+1,\ldots,n)}) = \bigcup\limits_{\Omega~\text{has even weight}}\overline{P_{\pi}^\Omega}$;\\
{\rm(2)}~$Supp^-_{T_\lambda}(f^{P_\pi}_{(n-m+1,\ldots,n)}) = \bigcup\limits_{\Omega~\text{has odd weight}}\overline{P_{\pi}^\Omega}.$
\end{lemma}
\proof
It follows from the definition of $PI$-eigenfunctions. $\square$

\medskip
The main theorem shows that each eigenfunction given by a polytabloid can be expressed as a sum of $PI$-eigenfunctions of $S_n^{JM}$.

\begin{theorem}\label{Thm2}
For a tableau $t$, the equality
\begin{equation}\label{maineq}
f_{\phi(\mathbf{e}_t)} = \sum\limits_{\sigma \in R_{t}(2) \times \ldots \times R_{t}(s)~}\sum\limits_{\pi \in CA_{\sigma(t)}} f^{P_\pi}_{(n-m+1,\ldots,n)}
\end{equation}
holds.
\end{theorem}
\proof
By Lemma \ref{SuppPolytabloid}, the function $f_{ \phi( \mathbf{e}_{t} ) }$ is a $(0,-1,1)$-function whose support can be partitioned into
the cosets. To prove that two $(0,-1,1)$-functions coincide, it is enough to show that their positive supports and their negative supports coincide. By Lemma \ref{SupPICoset}, the support of a $PI$-eigenfunction $f^{P_\pi}_{(n-m+1,\ldots,n)}$ can be partitioned into the cosets. By Lemma \ref{XandYCosets}, the codes of $2^{k}$ of them have no elements from $Y_t$ and
the codes of $2^m-2^{k}$ of them have at least one element from $Y_t$ (such cosets do not occur in the support of $f_{\phi(\mathbf{e}_t)}$).

Now we introduce a function $f_\sigma$, which represents the inner sum in the right part of the equality (\ref{maineq}), prove that it is a $(0,-1,1)$-function and investigate its support through cosets. For any $\sigma \in R_{t}(2) \times \ldots \times R_{t}(s)$, put
\begin{equation}\label{fsigma}
f_{\sigma} = \sum\limits_{\pi \in CA_{\sigma(t)}} f^{P_\pi}_{(n-m+1,\ldots,n)}.
\end{equation}
The following lemma gives an explicit partition of $Supp^+_{T_\lambda}(f_{\sigma})$ and $Supp^-_{T_\lambda}(f_{\sigma})$ into cosets.
\begin{lemma}\label{SuppRight}
For any $\sigma \in R_{t}(2) \times \ldots \times R_{t}(s)$, the function $f_\sigma$ is $(0,-1,1)$-function, where:
$$Supp^+_{T_\lambda}(f_{\sigma}) = \bigcup\limits_{\substack{\pi' \in C_{\sigma(t)},\\\pi'~\text{is even}} }
 [\pi'(t)];$$
$$Supp^-_{T_\lambda}(f_{\sigma}) = \bigcup\limits_{\substack{\pi' \in C_{\sigma(t)},\\\pi'~\text{is odd}} }
 [\pi'(t)].$$
\end{lemma}
\proof
It follows from Lemma~\ref{SupPICoset} that each of the sets $Supp^+_{T_\lambda}(f_{\sigma})$ and $Supp^-_{T_\lambda}(f_{\sigma})$
can be partitioned into the cosets such that, for each coset, $f_{\sigma}$ has the same value on permutations of this coset.
In fact, the supports of any two summands in the definition of $f_{\sigma}$ are disjoint.

From (\ref{fsigma}) we have that the support of $f_\sigma$ is included to the set
$$\bigcup\limits_{\substack{\pi \in CA_{\sigma(t)}, \\ \Omega \in \{0,1\}^m} }\overline{P_\pi^\Omega}.$$
Since $f_\sigma$ has the same value on permutations of a coset, let us investigate what is the value of $f_\sigma$ on a coset $H \in \{\overline{P_\pi^\Omega}~|~\pi \in CA_{\sigma(t)}, \Omega \in \{0,1\}^m\}$. Due to Lemma \ref{Cosets}, we have two equivalent ways to represent the coset $H$. Let us consider $H$ as the set of tableaux of shape $\lambda$ as follows:
$$
\ytableausetup{centertableaux, boxsize=1.7em}
\begin{ytableau}
 \ast & \ast & \ldots  & \ldots  & \ldots  & \ldots & \ldots & \ldots & \ldots & \ldots & \ast    \\
z_{12} & z_{22}  & \ldots & \ldots & \ldots & \ldots & \ldots & \ldots & z_{k2} \\
z_{13} & z_{23}  & \ldots & \ldots & \ldots & \ldots & z_{\lambda_33} \\
\ldots & \ldots & \ldots & \ldots & \ldots \\
z_{1s} & z_{2s} & \ldots & z_{\lambda_ss} \\
\end{ytableau}.
$$
In order to find the value of $f_\sigma$ on $H$ we have to determine all possible vectors $\Omega \in \{0,1\}^m$ and permutations $\pi \in CA_{\sigma(t)}$ such that the equality
\begin{equation}\label{CosetRestriction}
\overline{P_\pi^\Omega} = H
\end{equation}
holds. By Lemma \ref{CosetsAreEqual}, the cosets coincide if their codes are equal as vectors. This gives restrictions on possible values of $\Omega\in \{0,1\}^m$ and $\pi \in CA_{\sigma(t)}$. Since $\pi  = \pi_1 \circ \ldots \circ \pi_k$, we specify these restrictions columnwise, namely, for any $j \in \{1,\ldots,k\}$, we take the elements $z_{j2}, z_{j3}, \ldots, z_{j\mu_j}$ that correspond to $j$th column and occur in all tableaux from $H$, and obtain the system of $\mu_j-1$ equations on the permutation $\pi_j \in CA_{\sigma(t)}$:

\begin{equation}\label{Column}
  \begin{array}{ll}
    (\pi_{j1},\pi_{j2})^{\varepsilon_j} = z_{j2}, \\
    (y_{j3},\pi_{j3})^{\omega_{j3}} = z_{j3}, \\
    \multicolumn{1}{c}{\ldots}\\
    (y_{j\mu_j},\pi_{j\mu_j})^{\omega_{j\mu_j}} = z_{j\mu_j}.
  \end{array}
\end{equation}
Note that in the case $\mu_j = 2$ we have only one equation in (\ref{Column}).
This system implies that, for any $i \in \{3, \ldots, \mu_j\}$, we have
\begin{equation}\label{Rows3s}
\left\{
  \begin{array}{ll}
    \pi_{ji} = z_{ji}~\hbox{and}~\omega_{ji} = 0, & \hbox{if}~ z_{ji} \in X_{t}; \\
    y_{ji} = z_{ji}~\hbox{and}~\omega_{ji} = 1, & \hbox{if}~ z_{ji} \in Y_{t},
  \end{array}
\right.
\end{equation}
which follows from the fact that $\pi_{ji} \in X_t$ and
$y_{ji} \in Y_t$. Additionally, we have
\begin{equation}\label{Row2}
z_{j2} \in \{\pi_{j1}, \pi_{j2}\}.
\end{equation}
Let $CA_{\sigma(t)}^H(j)$ be the set of all permutations in $CA_{\sigma(t)}(j)$ satisfying the conditions~(\ref{Rows3s}) and~(\ref{Row2}). Note that $|CA_{\sigma(t)}^H(j)|$ is equal to the number of ways to assign $\pi_{ji}$ for all $i \in \{3, \ldots, \mu_j\}$ such that $z_{ji} \in Y_{t}$. In fact, if an even permutation has all but two already determined values, then it can be reconstructed uniquely. Thus, we have
\begin{equation}\label{AllAppropriate}
|CA_{\sigma(t)}^H(j)| = (|Z_j \cap Y_t|+1)!,
\end{equation}
 where $Z_j = \{ z_{j2}, \ldots, z_{j\mu_j}\}$.

Now our goal is to investigate the set $$CA_{\sigma(t)}^H = CA_{\sigma(t)}^H(1) \circ \ldots \circ CA_{\sigma(t)}^H(k),$$ which consists of all permutations $\pi$ from $CA_{\sigma(t)}$ whose support includes the coset $H$.

We consider the following two cases: $|Z_j \cap Y_t| > 0$ and $|Z_j \cap Y_t| = 0$.

\noindent
{\bf Case 1.}
Suppose that $|Z_j \cap Y_t| > 0$ holds, i.e. there exists $i \in \{3, \ldots, \mu_j\}$ such that $z_{ji} \in Y_{t}$.
Put $$(CA_{\sigma(t)}^H(j))' = \{\pi_j \in CA_{\sigma(t)}~|~\pi_{j1} = z_{j2}\}$$ and
$$(CA_{\sigma(t)}^H(j))'' = \{\pi_j \in CA_{\sigma(t)}~|~\pi_{j2} = z_{j2}\},$$ then
$$CA_{\sigma(t)}^H(j) = (CA_{\sigma(t)}^H(j))' \cup (CA_{\sigma(t)}^H(j))''.$$
Let us show that
\begin{equation}\label{HalfsColumn}
|(CA_{\sigma(t)}^H(j))'|  = |(CA_{\sigma(t)}^H(j))''|
\end{equation}
holds.
In fact, the two inclusions $$(CA_{\sigma(t)}^H(j))' \circ (\pi_{j1}~\pi_{j2}~\pi_{ji}) \subseteq (CA_{\sigma(t)}^H(j))''$$ and
$$(CA_{\sigma(t)}^H(j))'' \circ (\pi_{j2}~\pi_{j1}~\pi_{ji}) \subseteq (CA_{\sigma(t)}^H(j))'$$
hold, which gives (\ref{HalfsColumn}) since the equalities
$$|(CA_{\sigma(t)}^H(j))' \circ (\pi_{j1}~\pi_{j2}~\pi_{ji})| =
|(CA_{\sigma(t)}^H(j))'|$$ and
$$|(CA_{\sigma(t)}^H(j))'' \circ (\pi_{j2}~\pi_{j1}~\pi_{ji})| =
|(CA_{\sigma(t)}^H(j))''|$$
hold.
The equality (\ref{HalfsColumn}) is equivalent to the fact that, in (\ref{Column}), for one half of permutations in $CA_{\sigma(t)}^H(j)$, we have $\varepsilon_j = 1$ and, for another half, we have $\varepsilon_j = 0$.

\noindent
{\bf Case 2.}
Suppose that $|Z_j \cap Y_t| = 0$ holds. Then we have $|CA_{\sigma(t)}^H(j)| = 0$, and the permutation $\pi_j$
is uniquely determined. Hence, the value of $\varepsilon_j$ is uniquely determined in this case.

We have already proved that each factor $CA_{\sigma(t)}^H(j)$ in~(\ref{AllAppropriate}) either consists of one permutation
(in particular, $\varepsilon_j$ is uniquely determined) or, for a half of permutations in $CA_{\sigma(t)}^H(j)$, we have $\varepsilon_j = 1$ and, for another half, we have $\varepsilon_j = 0$. This implies that, if the code of $H$ has at least one element from $Y_t$, then for one half of permutations from $CA_{\sigma(t)}^H$ the vector $\Omega$ determined by~(\ref{CosetRestriction}) has odd number of ones, and for another half of permutations from $CA_{\sigma(t)}^H$
the vector $\Omega$ determined by~(\ref{CosetRestriction}) has even number of ones. Thus, $f_\sigma$ has zero value on the permutations of any coset $H$, whose code has at least one element from $Y_t$.

Finally, let us consider the case, when $H \in \{\overline{P_\pi^\Omega}~|~\pi \in CA_{\sigma(t)}\}$ is a coset, whose code has no elements from $Y_t$. We have proved that there exists a unique function from $\{f^{P_\pi}_{(n-m+1,\ldots,m)}~|~\pi \in CA_{\sigma(t)}\}$ that has non-zero value on the permutations of $H$. So, the function $f_\sigma$ has the same non-zero value ($1$ or $-1$) on $H$, and $H \subseteq Supp^+_{T_\lambda}(f_{\sigma}) \cup Supp^-_{T_\lambda}(f_{\sigma})$ holds.

To complete the proof of the lemma, it is enough to show that, for any $\pi' = \pi'_1 \circ \ldots \circ \pi'_k \in C_{\sigma(t)}$,
the inclusion $[\pi'(t)] \subseteq  Supp^+_{T_\lambda}(f_{\sigma})$ holds iff $\pi'$ is even. Put $H = [\pi'(t)]$ and take any $j \in \{1, \ldots, k\}$. Since $\pi'_{ji} \in X_t$, the condition (\ref{Rows3s}) gives $\pi_{ji} = \pi'_{ji}$ and, consequently, $\omega_j = 0$ for all $i \in \{3, \ldots, \mu_j\}$. The condition (\ref{Row2}) gives $\pi'_{j2} \in \{\pi_{j1}, \pi_{j2} \}$. Since $\pi_j$ is even, the equality $\pi_{j2} = \pi'_{j1}$ and, consequently, $\varepsilon_j = 1$ hold if $\pi'_j$ is odd. Thus, the number of ones in $\Omega$ is equal to the number of odd permutations $\pi'_j$, where $j$ runs over $\{1, \ldots, k\}$. This implies that $\pi'$ is even iff $\Omega$ has even number of ones, which proves the lemma. $\square$

By Lemmas~\ref{SuppPolytabloid} and~\ref{SuppRight}, the functions $f_{\phi(\mathbf{e}_t)}$ and $f_\sigma$ have the same values on the support of $f_\sigma$. For any distinct $\sigma_1, \sigma_2 \in R_{t}(2) \times \ldots \times R_{t}(s)$, the functions $f_{\sigma_1}$ and $f_{\sigma_2}$ have disjoint supports. The union of the supports of all functions $f_\sigma$, when $\sigma$ runs through $R_{t}(2) \times \ldots \times R_{t}(s)$, gives the support of $f_{\phi(\mathbf{e}_t)}$. The theorem is proved.
$\square$

\section{$PI$-eigenfunctions with eigenvalue $n-2$}\label{}
In this section, we prove that any eigenfunction corresponding to the largest non-principal eigenvalue $n-2$, can be reconstructed by its values on the second neighbourhood of a vertex.

Proposition~\ref{eigenfunctions} gives us the family of $PI$-eigenfunctions with eigenvalue $n-2$ when $m=1$. In this case we simplify~(\ref{MainFunction}) such that:

\begin{equation}\label{SimplyMainFunction}
f_{i}^{j,k}(\pi)=
\left\{
  \begin{array}{ll}
    1, & \hbox{if $\pi_j = i$;} \\
    -1, & \hbox{if $\pi_k = i$;}\\
     0, & \hbox{otherwise.}
  \end{array}
\right.
\end{equation}

Put
\begin{equation}\label{F2}
\mathcal{F}_2=\{f_i^{2,k}~|~i\in\{2,3,\ldots,n\}, k\in\{3,4,\ldots,n\}\}.
\end{equation}
Then the following statement holds.
\begin{lemma}\label{Basis} For $n \geqslant 3$, the set $\mathcal{F}_2$ forms a basis of the eigenspace of $S_n$ with eigenvalue $n-2$.
\end{lemma}
\proof For any $i\in\{2,3,\ldots,n\}$, the $PI$-eigenfunctions $f_i^{2,3}, f_i^{2,4}, \ldots, f_i^{2,n}$ are linearly independent since, for any $k\in\{3,4,\ldots,n\}$, the $PI$-eigenfunction $f_i^{2,k}$ has non-zero values on the set defined by $\pi_k = i$ only. It follows from the definition that, for any distinct $i_1,i_2 \in \{2,3,\ldots,n\}$ and any $k_1,k_2\in\{3,4,\ldots,n\}$, the $PI$-eigenfunctions $f_{i_1}^{2,k_1}$ and $f_{i_2}^{2,k_2}$ are orthogonal with respect to a natural inner product. This completes the proof. $\square$

Let us introduce additional notation. The second neighbourhood of the identity permutation in $S_n$ is presented by the set $N_2=\{(1rs)~|~r,s\in\{2,\ldots, n\}, r\ne s\},$ where $|N_2|=(n-1)(n-2)$ is the multiplicity of eigenvalue $n-2$~\cite{AKK16}.

We define a matrix $M_n$ as follows:
\begin{itemize}
  \item rows are indexed by elements from $\mathcal{F}_2$ divided into the following $n-1$ cohorts of length $n-2$ each:
$$f^{2,3}_2,f^{2,4}_2, f^{2,5}_2,\ldots, f^{2,n}_2~|~f^{2,3}_3,f^{2,4}_3, f^{2,5}_3,\ldots, f^{2,n}_3~|~
f^{2,4}_4,f^{2,3}_4, f^{2,5}_4,\ldots, f^{2,n}_4~|~$$ $$\ldots~|~\ldots~|~\ldots$$
$$~|~f^{2,i}_i,f^{2,3}_i, f^{2,4}_i,\ldots, f^{2,i-1}_i,f^{2,i+1}_i \ldots, f^{2,n}_{i}~|~\ldots~|~f^{2,n}_n,f^{2,3}_n, f^{2,5}_n,\ldots, f^{2,n}_{n-1}$$

  \item columns are indexed by the sequence of permutations from $N_2$ presented by the following $n-1$ blocks
of length $n-2$:
$$(132), (142), (152), \ldots, (1n2) ~|~(123), (143), (153), \ldots, (1n3)~|~$$
$$(124), (134), (154), \ldots, (1n4)~|~\ldots~|~(12n), (13n), (14n),\ldots,(1n-1n) $$
  \item the entries are the values of the row $PI$-eigenfunctions
on the corresponding column permutations.
\end{itemize}
In other words, the rows of the matrix $M_n$ are the restrictions of the $PI$-eigenfunctions from $\mathcal{F}_2$ on the set $N_2$. We consider the matrix $M_n$ as block-matrix $(B_{i_1,i_2})$, where $i_1,i_2 \in \{2,\ldots, n\}$, written in the form
$$M_n=
\left(
  \begin{array}{cccc}
    B_{2,2} & B_{2,3} & \ldots & B_{2,n-1} \\
    B_{3,2} & B_{3,3} & \ldots & B_{3,n-1} \\
    \vdots & \vdots & \ddots & \vdots \\
    B_{n,2} & B_{n,3} & \ldots & B_{n,n} \\
  \end{array}
\right).$$

The following  lemma can be proved by direct calculations.

\begin{lemma}\label{blocks} The blocks of the matrix $M_n$ are presented as follows.
\begin{enumerate}
  \item $B_{2,2} = -I_{n-2}$.
  \item $B_{2,i_2} = E_{n-2},$ where
  $ E_{n-2}=
  \left(
    \begin{array}{cccc}
      0 & 1 & \ldots & 1 \\
      \vdots & \vdots & ~ & \vdots \\
      0 & 1 & \ldots & 1 \\
    \end{array}
  \right)_{(n-2)\times(n-2)}
  $\\ for any  $i_2 \in \{3, \ldots, n\}$.

  \medskip
  \item $B_{i_1,i_2} = D_{n-2}$, where
  $D_{n-2} =
   \left(
     \begin{array}{ccccc}
       1 & 0 & 0 & \ldots & 0 \\
       1 & -1 & 0 & \ldots & 0 \\
       1 & 0 & -1 & \ldots & 0 \\
       \vdots & \vdots & \vdots & \ddots & 0 \\
       1 & 0 & 0 & \ldots & -1 \\
     \end{array}
   \right)_{(n-2)\times(n-2)}
  $\\ for any $i_1\in\{3, \ldots, n\}$, and $i_2 = i_1$.

 \medskip
  \item $B_{i_1,i_2} =
  \left\{
    \begin{array}{ll}
      C^{i_1-2}_{n-2}, & \hbox{if $i_1 > i_2$;} \\
      C^{i_1-1}_{n-2}, & \hbox{if $i_1 < i_2$,}
    \end{array}
  \right.$

  \medskip
   for any $i_1\in\{3, \ldots, n\}$, $i_2\in\{2, \ldots, n\}$, and $i_1 \ne i_2$,

  \medskip
  where $C_{n-2}^l =
    \left(
      \begin{array}{ccccccc}
        -1 & \ldots & -1 & 0 & -1 & \ldots & -1 \\
      0 & \ldots & 0 & 0 & 0 & \ldots & 0 \\
      \vdots & \vdots & \vdots & \vdots & \vdots & \vdots & \vdots \\
      0 & \ldots & 0 & 0 & 0 & \ldots & 0 \\

      \end{array}
    \right)_{(n-2)\times(n-2)}
  $
  with zero $l$-th column, $l \in \{1, \ldots, n-2\}$.
\end{enumerate}
\end{lemma}
$\square$

\medskip

By Lemma~\ref{blocks}, we have
$$
M_n=
\left(
  \begin{array}{cccccc}
    -I_{n-2} & E_{n-2} & E_{n-2}  & \ldots & E_{n-2} & E_{n-2} \\
    C^{1}_{n-2} & D_{n-2} & C^{2}_{n-2} & \ldots & C^{2}_{n-2} & C^{2}_{n-2} \\
    C^{2}_{n-2} & C^{2}_{n-2} & D_{n-2} & \ldots & C^{3}_{n-2} & C^{3}_{n-2}\\
    \vdots & \vdots & \vdots & \ddots & \vdots  & \vdots \\
    C^{n-3}_{n-2} & C^{n-3}_{n-2} & C^{n-3}_{n-2}  & \ldots & D_{n-2} & C^{n-2}_{n-2} \\
    C^{n-2}_{n-2} & C^{n-2}_{n-2} & C^{n-2}_{n-2}  & \ldots & C^{n-2}_{n-2} & D_{n-2} \\
  \end{array}
\right)
$$

The following lemma gives us an explicit  formula for the determinant $det(M_n)$ of the matrix $M_n$. A column of blocks (row of blocks) in the matrix $M_n$ is called a {\it block-column} (a {\it block-row}).

\begin{lemma}\label{det}
The following equality holds
\begin{equation}\label{e7}
det(M_n)=(-1)^{n-3}(n-2)^{n-2}(n^2-5n+5).
\end{equation}
In particular, the matrix $M_n$ is non-degenerate.
\end{lemma}
\proof We prove~(\ref{e7}) applying step-by-step operations to the matrix $M_n$ with $n-1$ block-columns of size $n-2$. Those matrix transformations are presented in Appendix. $\square$

\medskip
Since the rows of the matrix $M_n$ are the restrictions of the $PI$-eigenfunctions from $\mathcal{F}_2$ on the set $N_2$, where $|N_2|=|\mathcal{F}_2|=(n-1)(n-2)$, from Lemma~\ref{Basis}, Lemma~\ref{det} and the vertex-transitivity of $S_n$ we have the following result.

\begin{theorem} Any eigenfunction of $S_n$, $n \geqslant 3$, with eigenvalue $n-2$ can be uniquely reconstructed by its values on the second neighbourhood of any vertex.
\end{theorem}

The reconstruction problem for eigenfunctions and perfect codes of the Hamming graph was investigated in~\cite{A95,AV08,H02,V15}.

\section*{Acknowledgment} \label{sec3}

The authors are grateful to Sergey Avgustinovich and Akihiro Munemasa for fruitful discussions in Japan, China and Russia. Their remarkable questions and  useful suggestions significantly improved the content of this paper. Special thanks goes to Alexander A.~Ivanov for careful reading of the manuscript at its final stage.

\setlength{\tabcolsep}{1.8pt}

\begin{table}
\begin{flushright}
APPENDIX
\end{flushright}
\begin{turn}{90}

\end{turn}
\begin{turn}{90}
Step 1: Take the first column of the second block-column, add all other columns to it, factor out $n^2-5n+5$.
\end{turn}
\end{table}

\begin{table}
\begin{flushright}
APPENDIX
\end{flushright}
\begin{turn}{90}
\end{turn}
\begin{turn}{90}
Step 2: Take the last column of the first block-column, add other columns of this block-column and the first column

\end{turn}
\begin{turn}{90}
 of the second block-column to it, factor out $n-2$.
\end{turn}
\end{table}

\begin{table}
\begin{flushright}
APPENDIX
\end{flushright}
\begin{turn}{90}
\end{turn}
\begin{turn}{90}
 Step 3: Take the last column of the first block-column, subtract it from the first column of the second block-column.
\end{turn}
\end{table}

\begin{table}
\begin{flushright}
APPENDIX
\end{flushright}
\begin{turn}{90}
\end{turn}
\begin{turn}{90}
Step 4: For all but the first columns of all but the first block-columns,

\end{turn}
\begin{turn}{90}
subtract the first column of the second block-column from them.
\end{turn}
\end{table}

\begin{table}
\begin{flushright}
APPENDIX
\end{flushright}
\begin{turn}{90}
\begin{tabular}{ccccc>{\columncolor{blue!25}}c|>{\columncolor{blue!25}}cccccc|cccccc|cccccc|cccccc|cccccc|cccccc}

-1& 0 & 0 & \dots & 0 & 0  &  1 & 0 & 0 & \dots & 0 & 0 & 0 & 0 & 0 & \dots & 0 & 0 & 0 & 0 & 0 & \dots & 0 & 0 & & & &  & & &   0 & 0 & 0 & \dots & 0 & 0 &     0 & 0 & 0 & \dots & 0 & 0  \\
0 & -1 & 0 & \dots & 0 & 0  & 1 & 0 & 0 & \dots & 0 & 0 & 0 & 0 & 0 & \dots & 0 & 0 & 0 & 0 & 0 & \dots & 0 & 0 &  &  &  &  &  &  & 0 & 0 & 0 & \dots & 0 & 0 & 0 & 0 & 0 & \dots & 0 & 0  \\
$\cdot$ & $\cdot$ & $\cdot$ &  & $\cdot$ & $\cdot$  & $\cdot$ & $\cdot$ & $\cdot$ &  & $\cdot$ & $\cdot$ & $\cdot$ & $\cdot$ & $\cdot$ &  & $\cdot$ & $\cdot$ & $\cdot$ & $\cdot$ & $\cdot$ &  & $\cdot$ & $\cdot$ &  &  & & $\dots$ &  & & $\cdot$ & $\cdot$ & $\cdot$ &  & $\cdot$ & $\cdot$ & $\cdot$ & $\cdot$ & $\cdot$ &  & $\cdot$ & $\cdot$  \\
0 & 0 & 0 & \dots & 0 & 0  & 1 & 0 & 0 & \dots & 0 & 0 & 0 & 0 & 0 & \dots & 0 & 0 & 0 & 0 & 0 & \dots & 0 & 0 &  &  &  &  &  &  & 0 & 0 & 0 & \dots & 0 & 0 & 0 & 0 & 0 & \dots & 0 & 0  \\
\hline
0& -1 & -1 & \dots & -1 & -1  & 0 & 0 & 0 & \dots & 0 & 0 & -1 & 0 & -1 & \dots & -1 & -1 & -1 & 0 & -1 & \dots & -1 & -1 & & & &  & & &   -1 & 0 & -1 & \dots & -1 & -1 &     -1 & 0 & -1 & \dots & -1 & -1  \\
0 & 0 & 0 & \dots & 0 & 0  & 0 & -1 & 0 & \dots & 0 & 0 & 0 & 0 & 0 & \dots & 0 & 0 & 0 & 0 & 0 & \dots & 0 & 0 &  &  &  &  &  &  & 0 & 0 & 0 & \dots & 0 & 0 & 0 & 0 & 0 & \dots & 0 & 0  \\
0 & 0 & 0 & \dots & 0 & 0  & 0 & 0 & -1 & \dots & 0 & 0 & 0 & 0 & 0 & \dots & 0 & 0 & 0 & 0 & 0 & \dots & 0 & 0 &  &  &  &  &  &  & 0 & 0 & 0 & \dots & 0 & 0 & 0 & 0 & 0 & \dots & 0 & 0  \\
$\cdot$ & $\cdot$ & $\cdot$ &  & $\cdot$ & $\cdot$  & $\cdot$ & $\cdot$ & $\cdot$ &  & $\cdot$ & $\cdot$ & $\cdot$ & $\cdot$ & $\cdot$ &  & $\cdot$ & $\cdot$ & $\cdot$ & $\cdot$ & $\cdot$ &  & $\cdot$ & $\cdot$ &  &  & & $\dots$ &  & & $\cdot$ & $\cdot$ & $\cdot$ &  & $\cdot$ & $\cdot$ & $\cdot$ & $\cdot$ & $\cdot$ &  & $\cdot$ & $\cdot$  \\
\end{tabular}
\end{turn}
\begin{turn}{90}
Step 5: Interchange the last column of the first block-column with the first column of the
second block-column,
\end{turn}
\begin{turn}{90}
factor out $-1$.
\end{turn}
\end{table}

\begin{table}
\begin{flushright}
APPENDIX
\end{flushright}
\begin{turn}{90}
\begin{tabular}{cccccc|cccccc|>{\columncolor{blue!25}}cccccc|>{\columncolor{blue!25}}cccccc|cccccc|>{\columncolor{blue!25}}cccccc|>{\columncolor{blue!25}}cccccc}

-1& 0 & 0 & \dots & 0 & 1  &  0 & 0 & 0 & \dots & 0 & 0 & 0 & 0 & 0 & \dots & 0 & 0 & 0 & 0 & 0 & \dots & 0 & 0 & & & &  & & &   0 & 0 & 0 & \dots & 0 & 0 &     0 & 0 & 0 & \dots & 0 & 0  \\
0 & -1 & 0 & \dots & 0 & 1  & 0 & 0 & 0 & \dots & 0 & 0 & 0 & 0 & 0 & \dots & 0 & 0 & 0 & 0 & 0 & \dots & 0 & 0 &  &  &  &  &  &  & 0 & 0 & 0 & \dots & 0 & 0 & 0 & 0 & 0 & \dots & 0 & 0  \\
$\cdot$ & $\cdot$ & $\cdot$ &  & $\cdot$ & $\cdot$  & $\cdot$ & $\cdot$ & $\cdot$ &  & $\cdot$ & $\cdot$ & $\cdot$ & $\cdot$ & $\cdot$ &  & $\cdot$ & $\cdot$ & $\cdot$ & $\cdot$ & $\cdot$ &  & $\cdot$ & $\cdot$ &  &  & & $\dots$ &  & & $\cdot$ & $\cdot$ & $\cdot$ &  & $\cdot$ & $\cdot$ & $\cdot$ & $\cdot$ & $\cdot$ &  & $\cdot$ & $\cdot$  \\
0 & 0 & 0 & \dots & 0 & 1  & 0 & 0 & 0 & \dots & 0 & 0 & 0 & 0 & 0 & \dots & 0 & 0 & 0 & 0 & 0 & \dots & 0 & 0 &  &  &  &  &  &  & 0 & 0 & 0 & \dots & 0 & 0 & 0 & 0 & 0 & \dots & 0 & 0  \\
\hline
0& -1 & -1 & \dots & -1 & -0  & -1 & 0 & 0 & \dots & 0 & 0 & -1 & 0 & -1 & \dots & -1 & -1 & -1 & 0 & -1 & \dots & -1 & -1 & & & &  & & &   -1 & 0 & -1 & \dots & -1 & -1 &     -1 & 0 & -1 & \dots & -1 & -1  \\
0 & 0 & 0 & \dots & 0 & 0  & 0 & -1 & 0 & \dots & 0 & 0 & 0 & 0 & 0 & \dots & 0 & 0 & 0 & 0 & 0 & \dots & 0 & 0 &  &  &  &  &  &  & 0 & 0 & 0 & \dots & 0 & 0 & 0 & 0 & 0 & \dots & 0 & 0  \\
0 & 0 & 0 & \dots & 0 & 0  & 0 & 0 & -1 & \dots & 0 & 0 & 0 & 0 & 0 & \dots & 0 & 0 & 0 & 0 & 0 & \dots & 0 & 0 &  &  &  &  &  &  & 0 & 0 & 0 & \dots & 0 & 0 & 0 & 0 & 0 & \dots & 0 & 0  \\
\end{tabular}
\end{turn}
\begin{turn}{90}
Step 6: For each block-column excepting the first and the second ones,
\end{turn}
\begin{turn}{90}
add non-first columns of this block-column to its first column.
\end{turn}
\end{table}

\begin{table}
\begin{flushright}
APPENDIX
\end{flushright}
\begin{turn}{90}
\begin{tabular}{cccccc|cccccc|>{\columncolor{blue!25}}cccccc|>{\columncolor{blue!25}}cccccc|cccccc|>{\columncolor{blue!25}}cccccc|>{\columncolor{blue!25}}cccccc}

-1& 0 & 0 & \dots & 0 & 1  &  0 & 0 & 0 & \dots & 0 & 0 & 0 & 0 & 0 & \dots & 0 & 0 & 0 & 0 & 0 & \dots & 0 & 0 & & & &  & & &   0 & 0 & 0 & \dots & 0 & 0 &     0 & 0 & 0 & \dots & 0 & 0  \\
0 & -1 & 0 & \dots & 0 & 1  & 0 & 0 & 0 & \dots & 0 & 0 & 0 & 0 & 0 & \dots & 0 & 0 & 0 & 0 & 0 & \dots & 0 & 0 &  &  &  &  &  &  & 0 & 0 & 0 & \dots & 0 & 0 & 0 & 0 & 0 & \dots & 0 & 0  \\
$\cdot$ & $\cdot$ & $\cdot$ &  & $\cdot$ & $\cdot$  & $\cdot$ & $\cdot$ & $\cdot$ &  & $\cdot$ & $\cdot$ & $\cdot$ & $\cdot$ & $\cdot$ &  & $\cdot$ & $\cdot$ & $\cdot$ & $\cdot$ & $\cdot$ &  & $\cdot$ & $\cdot$ &  &  & & $\dots$ &  & & $\cdot$ & $\cdot$ & $\cdot$ &  & $\cdot$ & $\cdot$ & $\cdot$ & $\cdot$ & $\cdot$ &  & $\cdot$ & $\cdot$  \\
0 & 0 & 0 & \dots & 0 & 1  & 0 & 0 & 0 & \dots & 0 & 0 & 0 & 0 & 0 & \dots & 0 & 0 & 0 & 0 & 0 & \dots & 0 & 0 &  &  &  &  &  &  & 0 & 0 & 0 & \dots & 0 & 0 & 0 & 0 & 0 & \dots & 0 & 0  \\
\hline
0& -1 & -1 & \dots & -1 & -0  & -1 & 0 & 0 & \dots & 0 & 0 & -(n-3) & 0 & -1 & \dots & -1 & -1 & -(n-3) & 0 & -1 & \dots & -1 & -1 & & & &  & & &   -(n-3) & 0 & -1 & \dots & -1 & -1 &     -(n-3) & 0 & -1 & \dots & -1 & -1  \\
0 & 0 & 0 & \dots & 0 & 0  & 0 & -1 & 0 & \dots & 0 & 0 & 0 & 0 & 0 & \dots & 0 & 0 & 0 & 0 & 0 & \dots & 0 & 0 &  &  &  &  &  &  & 0 & 0 & 0 & \dots & 0 & 0 & 0 & 0 & 0 & \dots & 0 & 0  \\
0 & 0 & 0 & \dots & 0 & 0  & 0 & 0 & -1 & \dots & 0 & 0 & 0 & 0 & 0 & \dots & 0 & 0 & 0 & 0 & 0 & \dots & 0 & 0 &  &  &  &  &  &  & 0 & 0 & 0 & \dots & 0 & 0 & 0 & 0 & 0 & \dots & 0 & 0  \\
\end{tabular}
\end{turn}
\begin{turn}{90}
Step 7: For each block-column excepting the first and the second ones,

\end{turn}
\begin{turn}{90}
take the first column and add the first column of the second block-column multiplied by $-(n-3)$ to it;
\end{turn}
\begin{turn}{90}
factor out $n-2$ from each of the modified $n-3$ first columns.
\end{turn}
\end{table}

\begin{table}
\begin{flushright}
APPENDIX
\end{flushright}
\begin{turn}{90}
\begin{tabular}{cccccc|c>{\columncolor{blue!25}}c>{\columncolor{blue!25}}c>{\columncolor{blue!25}}c>{\columncolor{blue!25}}c>{\columncolor{blue!25}}c|c>{\columncolor{blue!25}}c>{\columncolor{blue!25}}c>{\columncolor{blue!25}}c>{\columncolor{blue!25}}c>{\columncolor{blue!25}}c|c>{\columncolor{blue!25}}c>{\columncolor{blue!25}}c>{\columncolor{blue!25}}c>{\columncolor{blue!25}}c>{\columncolor{blue!25}}c|cccccc|c>{\columncolor{blue!25}}c>{\columncolor{blue!25}}c>{\columncolor{blue!25}}c>{\columncolor{blue!25}}c>{\columncolor{blue!25}}c|c>{\columncolor{blue!25}}c>{\columncolor{blue!25}}c>{\columncolor{blue!25}}c>{\columncolor{blue!25}}c>{\columncolor{blue!25}}c}

-1& 0 & 0 & \dots & 0 & 1  &  0 & 0 & 0 & \dots & 0 & 0 & 0 & 0 & 0 & \dots & 0 & 0 & 0 & 0 & 0 & \dots & 0 & 0 & & & &  & & &   0 & 0 & 0 & \dots & 0 & 0 &     0 & 0 & 0 & \dots & 0 & 0  \\
0 & -1 & 0 & \dots & 0 & 1  & 0 & 0 & 0 & \dots & 0 & 0 & 0 & 0 & 0 & \dots & 0 & 0 & 0 & 0 & 0 & \dots & 0 & 0 &  &  &  &  &  &  & 0 & 0 & 0 & \dots & 0 & 0 & 0 & 0 & 0 & \dots & 0 & 0  \\
$\cdot$ & $\cdot$ & $\cdot$ &  & $\cdot$ & $\cdot$  & $\cdot$ & $\cdot$ & $\cdot$ &  & $\cdot$ & $\cdot$ & $\cdot$ & $\cdot$ & $\cdot$ &  & $\cdot$ & $\cdot$ & $\cdot$ & $\cdot$ & $\cdot$ &  & $\cdot$ & $\cdot$ &  &  & & $\dots$ &  & & $\cdot$ & $\cdot$ & $\cdot$ &  & $\cdot$ & $\cdot$ & $\cdot$ & $\cdot$ & $\cdot$ &  & $\cdot$ & $\cdot$  \\
0 & 0 & 0 & \dots & 0 & 1  & 0 & 0 & 0 & \dots & 0 & 0 & 0 & 0 & 0 & \dots & 0 & 0 & 0 & 0 & 0 & \dots & 0 & 0 &  &  &  &  &  &  & 0 & 0 & 0 & \dots & 0 & 0 & 0 & 0 & 0 & \dots & 0 & 0  \\
\hline
0& -1 & -1 & \dots & -1 & -0  & -1 & 0 & 0 & \dots & 0 & 0 & 0 & 0 & -1 & \dots & -1 & -1 & 0 & 0 & -1 & \dots & -1 & -1 & & & &  & & &   0 & 0 & -1 & \dots & -1 & -1 &     0 & 0 & -1 & \dots & -1 & -1  \\
\rowcolor{blue!25}
0 & 0 & 0 & \dots & 0 & 0  & 0 & -1 & 0 & \dots & 0 & 0 & 0 & 0 & 0 & \dots & 0 & 0 & 0 & 0 & 0 & \dots & 0 & 0 &  &  &  &  &  &  & 0 & 0 & 0 & \dots & 0 & 0 & 0 & 0 & 0 & \dots & 0 & 0  \\
\end{tabular}
\end{turn}
\begin{turn}{90}
Step 8: For all block-rows but the first one,
\end{turn}
\begin{turn}{90}
expand the determinant by all their rows but the first one;
\end{turn}
\begin{turn}{90}
obtain the product of $(-1)^{(n-2)(n-3)}=1$ and the minor.
\end{turn}
\end{table}

\newpage
\begin{flushright}
APPENDIX
\end{flushright}
\begin{tabular}{cccccc|c|c|c|c|c|c}

-1& 0 & 0 & \dots & 0 & 1  &  0 &  0 & 0  & &    0  &    0   \\
0 & -1 & 0 & \dots & 0 & 1  &  0 &  0 & 0  & &    0  &    0   \\
$\cdot$ & $\cdot$ & $\cdot$ &  & $\cdot$ & $\cdot$  & $\cdot$ & $\cdot$ & $\cdot$ & \dots & $\cdot$& $\cdot$\\
0 & 0 & 0 & \dots & -1 & 1  &  0 &  0 & 0  & &    0   &   0   \\
0 & 0 & 0 & \dots & 0 & 1  &  0 &  0 & 0  & &    0   &   0   \\
\hline
0 & -1 & -1 & \dots & -1 & 0  &  -1 &  0 & 0  & &    0   &   0   \\
\hline
-1 & 0 & -1 & \dots & -1 & 0  &  -1 &  1 & 0  & &    0   &   0   \\
\hline
$\cdot$ & $\cdot$ & $\cdot$ &  & $\cdot$ & $\cdot$  & $\cdot$ & $\cdot$ & $\cdot$ & & $\cdot$ & $\cdot$ \\
\hline
-1 & -1 & -1 & \dots & 0 & 0  &  -1 &  0 & 0  & &    1   &   0   \\
\hline
-1 & -1 & -1 & \dots & -1 & 0  &  -1 &  0 & 0  & &    0   &   1   \\

\end{tabular}
\\
Step 9: The reduced matrix is a block triangular matrix.


\section*{References}

\begin{thebibliography}{00}

\bibitem{A95}
S~.V.~Avgustinovich, On One Property of Perfect Binary Codes, Diskretn. Anal. Issled. Oper., \textbf{2} (1) (1995) 4--6.

\bibitem{AV08}
S.~V.~Avgustinovich, A.~Yu.~Vasil’eva, Reconstruction theorems for centered functions and perfect codes, Siberian Mathematical Journal, \textbf{49} (3) (2008)  383--388.

\bibitem{AKK16}
S.~V.~Avgustinovich, E.~N.~Khomyakova, E.~V.~Konstantinova, Multiplicities of eigenvalues of the Star graph,
\emph{Siberian Electronic Mathematical Report}, \textbf{13} (2016) 1258--1270.

\bibitem{BH12}
A.~E.~Brouwer, W.~H.~Haemers, \emph{Spectra of Graphs}, Springer, New York, 2012.

\bibitem{CF12}
G.~Chapuy, V.~Feray, A note on a Cayley graph of $\mathrm{Sym}_{n}$, \emph{arXiv:1202.4976v2} (2012) 1--3.

\bibitem{FOW85}
L.~Flatto, A.~Odlyzko, D.~Wales, Random shuffles and group representations. \emph{Ann. Prob.}, \textbf{13} (1) (1985)  154--178.

\bibitem{H02}
O.~Heden, On the Reconstruction of Perfect Codes, Discrete Mathematics, \textbf{256} (1--2) (2002) 479--485.

\bibitem{J74}
A.~Jucys, Symmetric polynomials and the center of the symmetric group ring, \emph{Reports Math. Phys.}, \textbf{5} (1974) 107--112.

\bibitem{KM75}
R.~Krakovski, B.~Mohar, Spectrum of Cayley graphs on the symmetric group generated by transposition, \emph{Linear Algebra and its Applications}, \textbf{437} (2012) 1033--1039.

\bibitem{KMP16}
D.~Krotov, I.~Mogilnykh, V.~Potapov, To the theory of q-ary Steiner and other-type trades, \emph{Discrete Mathematics}, \textbf{339} (3) (2016) 1150--1157.

\bibitem{K17}
D.~Krotov, The extended $1$-perfect trades in small hypercubes, \emph{Discrete Mathematics}, \textbf{340} (10) (2017) 2559--2572.

\bibitem{R11}
P.~Renteln, The distance spectra of Cayley graphs of Coxeter groups, \emph{Discrete Mathematics}, \textbf{311} (2011) 738--755.

\bibitem{Sa01}
B.~E.~Sagan, \emph{The Symmetric Group: representations, combinatorial algorithms, and symmetric functions}, Second ed. \emph{Graduate texts in mathematics.}, \textbf{203} (2001) 240~pp.

\bibitem{V17}
A.~Valyuzhenich, Minimum supports of eigenfunctions of Hamming graphs, \emph{Discrete Mathematics}, \textbf{340} (5) (2017) 1064--1068.

\bibitem{V15}
A.~Yu.~Vasil’eva, Reconstruction of eigenfunctions of a $q$-ary $n$-dimensional hypercube, Problems of Information Transmission, \textbf{51} (3) (2015) 231--239.

\end{thebibliography}
\end{document}